%% file: 523.tex
\ifx\shlhetal\undefinedcontrolsequence\let\shlhetal\relax\fi


\input amstex
\input mathdefs
\sectno=-1   
\localtags
\NoBlackBoxes
\def\cite #1{\rm[#1]} 
\documentstyle {amsppt}
\topmatter
\title {{\it Almost Free groups Below the First Fix Point Exist}\\
Sh523} \endtitle
\rightheadtext {Almost Free Groups}
\author {Saharon Shelah \thanks {\null\newline
   Research supported by
   German-Israeli Foundation for Scientific Research \& Development
   Grant No. G-294.081.06/93}
\endthanks } \endauthor
\affil {Institute of Mathematics \\
The Hebrew University\\
Jerusalem, Israel
\medskip
Rutgers University\\
Department of Mathematics\\
New Brunswick, NJ USA} \endaffil
\endtopmatter
\medskip
\document

\specialhead {Anotated Content} \endspecialhead
\resetall
\medskip
\noindent
\S1 - $\lambda$-free does not imply free for $\lambda <$ first fix point
\newline
[We represent material from \cite{MgSh:204}, assuming knowledge of
\cite{Sh:g, Ch.II}. \newline
Just enough to make \S2 intelligible.]
\medskip

\noindent
\S2 - Nicely compact sets and non-reflection \newline
[We answer a question of Foreman and Magidor on reflection of stationary
subsets of ${\Cal S}_{< \aleph_2}(\lambda) = \{ a \subseteq \lambda:
|a| < \aleph_2 \}$].  We answer a question of Mekler Eklof on the 
closure operations of the incompactness spectrum.]
\medskip

\noindent
\S3 - $NPT$ is not transitive \newline
[We prove the consistency of $NPT(\lambda,\mu) + NPT(\mu,\kappa)
\nRightarrow NPT(\lambda,\kappa)$].
\newpage

\head {\S1 $\lambda$-Free Does Not Imply Free for $\lambda <$ First
Fix Point} \endhead
\resetall
\bigskip

\definition {\stag{1.0} Definition}  $NPT(\lambda,\kappa)$ means that 
there is a family
$\{ A_i:i < \lambda \}$ of subsets of $\lambda,|A_i| \le \kappa$, the family
has no transversal but for $\alpha < \lambda$ we have $\{ A_i:i < \alpha \}$
has a transversal.  Let $PT(\lambda,\kappa)$ be $\neg NPT(\lambda,\kappa)$.
\newline
Now by \cite{Sh:161}, $NPT(\lambda,\aleph_0)$ is equivalent to the 
existence of $\lambda$-free not $\lambda^+$-free abelian group.  
Also ``canonical" examples of ``incompactness" are given there.
\enddefinition
\bigskip

\noindent
Magidor Shelah \cite{MgSh:204} proved
\proclaim{\stag{1.1} Theorem (ZFC)}  If $\lambda$ is smaller than the first
$\alpha = \aleph_\alpha > \aleph_0$ \underbar{then} there is a $\lambda$-free
not free abelian group.
\medskip
\noindent
This was done by providing an induction step from incompactness in
$\kappa$, \newline
$\aleph_0 < \kappa = \text{ cf}(\kappa) < \aleph_\kappa$ to 
incompactness in some $\lambda \ge \aleph_\kappa$, in fact 
$\aleph_{\kappa+1}$. \newline
However, subsequently, to answer questions of Eklof and Mekler 
(in \scite{2.2}, \scite{2.3}) and of Magidor and Foreman (in \scite{2.7}), 
we do more.
\medskip

\noindent 
\underbar{Explanation of the History}:  Generally on the history
of this problem see the book of Eklof and Mekler \cite{EM}.  The 
existence of such group for $\lambda = \aleph_{\omega+1}$ is proved in 
\cite{Sh:108} assuming $\aleph_\omega$ strong
limit or just $2^{\aleph_0} < \aleph_\omega$; this is based on investigating
$I[\aleph_{\omega+1}]$ (where for a regular uncountable $\lambda$, \newline
$I[\lambda] =
\{S \subseteq \lambda:\text{for some club } E \text{ of } \lambda \text{ and }
\bar a = \langle a_\alpha:\alpha < \lambda \rangle \text{ we have }
a_\alpha \subseteq \alpha,\text{ otp}(a_\alpha) < \alpha \text{ and }
\alpha \in S \cap E \Rightarrow \alpha = \sup(a_\alpha)\}$), see on it
\cite{Sh:g, {AG \S1}}.  In the earlier version of \cite{MgSh:204}, 
for $\lambda > \text{ cf}(\lambda) = \kappa$ when $2^\kappa < \lambda$ we 
replace a division to cases defined from properties $I[\lambda]$ by one 
defined from 
$\bar f = \langle f_\alpha:\alpha < \lambda^+ \rangle <_I$-increasing
cofinal in $\dsize \prod_{i < \kappa} \lambda_i/I$. \newline
It is shown there that the set of bad (for ${\bar f}$) ordinals
$< \aleph_{\kappa+1}$ is ``small", (but assuming something like 
$2^\kappa < \lambda$, then an example in $\lambda$ is gotten from one on 
$\kappa^+$.  In the present version of \cite{MgSh:204} this was eliminated
by showing: for good $\delta < \lambda$ there is a club of good ordinals 
$< \delta$ (for ${\bar f}$); and similarly for bad (we ignore ordinals of 
small cofinality).
\endproclaim
\bigskip

\noindent
In this section we represent analyses from \cite{MgSh:204} needed in \S2; 
it is organized around the proof of
\proclaim{\stag{1.2} Lemma}  If NPT$(\kappa,\theta)$, and $\kappa \ge \theta =
\text{ cf}(\theta)$ \underbar{then} NPT$(\kappa^{+(\kappa+1)},\theta)$.
\endproclaim
\bigskip

\remark{Remark}  1) Similarly NPT$_I$, see on it \cite{Sh:355, {\S6}} 
(for $I$ an ideal on $\theta$ extending $J^{bd}_\theta$).
\endremark
\bigskip

\noindent
Till end of the proof of 1.2 we make:
\demo{\stag{1.3} Convention}  $\kappa$ is a fixed regular cardinal.  
$\lambda$ is a cardinal $> \kappa$ of cofinality $\kappa$.

(For proving \scite{1.1}, we need just $\lambda = \aleph_\kappa$, but 
latter we will be interested in other values too).  $J$ denotes ideals 
on $\kappa$ extending $J^{bd}_\kappa = \{A \subseteq \kappa:A 
\text{ bounded}\}$.\newline
$D,E$ denotes filters on $\kappa$ containing the co-bounded subsets of
$\lambda$.
\medskip
The rest of the section proves \scite{1.2}.
\enddemo
\bigskip

\proclaim{\stag{1.4} Claim}  Suppose 
$\lambda = \kappa^{+\kappa}$, $\mu = \lambda^+;J = J^{bd}_\kappa$, 
\underbar{then} there are $\langle \lambda_i:i < \kappa \rangle$, \newline
$\langle f_\alpha:\alpha < \mu \rangle$ such that the following holds:
\medskip
\roster
%
\item "{$(*)(a)$}"  $\langle \lambda_i:i < \kappa \rangle \text{ strictly
increasing sequence of regular cardinals } > \kappa^+$ \newline
$\text{ with limit } \lambda$,
\smallskip
\noindent
\item "{$(b)$}"  for $\alpha < \mu,f_\alpha \in \Pi \lambda_i,
[\alpha < \beta < \mu \Rightarrow f_\alpha <_J f_\beta]$,
\smallskip
\noindent
\item "{$(c)$}"  $\left( \forall f \in \dsize \prod_{i<\kappa} \lambda_i 
\right) (\exists \alpha < \mu)[f <_J f_\alpha]$,
\smallskip
\noindent
\item "{$(d)$}"  $\mu = \text{ cf}(\mu) > \lambda$ (what is trivial here),
\smallskip
\noindent
\item "{$(e)$}"  $J \text{ is an ideal on } \kappa \text{ containing } 
J^{bd}_\kappa$.
\endroster
\endproclaim 
\bigskip

\demo{Proof}  By \cite{Sh:355, {1.5}}.
\enddemo
\bigskip

\centerline{$* \qquad * \qquad *$}
\bigskip

\noindent
For awhile (till we finish proving \scite{1.12}) we shall assume $(*)$ of 
\scite{1.4}.
\demo{\stag{1.5} Context}  $(*)$ of \scite{1.4}.
\enddemo
\bigskip

\noindent
Let us quote \cite{Sh:355, {1.2}}, for the reader's convenience.
\proclaim{\stag{1.6A}Claim}  Assume cf$(\delta) > \kappa^+,I$ is an ideal on
$\kappa$ and suppose $\langle f_\alpha:\alpha < \delta \rangle$ is a
$<_I$-increasing sequence of members of ${}^\kappa\text{Ord}$.
\underbar{Then} exactly one of the following holds:
\medskip
\roster
%
\item "{$(i)$}"  for some ultrafilter $D$ on $\kappa$ disjoint from $I$ we
have:
{\roster
\itemitem{ $(*)_D$ }  there are sets $s_i \subseteq \text{ Ord}, |s_i| \le
\kappa \text{ and } \langle \alpha_\zeta:\zeta < \text{ cf}(\delta) \rangle$
increasing continuous with limit $\delta$, such that for each $\zeta <
\text{ cf}(\delta)$ for some \newline
$h_\zeta \in \dsize \prod_{i < \kappa} s_i$ we have: \newline
$f_{\alpha_\zeta}/D < h_\zeta/D < f_{\alpha_{\zeta +1}}/D$
\endroster}
\item "{$(ii)$}"  $(**)_I \quad$ some $f \in {}^\kappa\text{Ord}$ is a
$<_I$-eub of $\langle f_\alpha:\alpha < \delta \rangle$, \newline

$\qquad \quad$ (i.e. $f$ satisfies $(\alpha) + (\beta)$ below) and
$(\gamma)$ holds:
{\roster
\itemitem{ $(\alpha)$ }  for $\alpha < \delta,f_\alpha <_I f$
\itemitem{ $(\beta)$ }  if $g \in {}^\kappa\text{Ord},g <_I f$ then for some
$\alpha,g <_I f_\alpha$ and
\itemitem{ $(\gamma)$ }  cf$[f(i)] > \kappa \text{ for } i < \kappa$
\endroster}
\item "{$(iii)$}"  condition $(i)$ fails and \newline
$(***)_I \quad$ for some unbounded $A \subseteq \delta$ and $t_\alpha
\subseteq \kappa$ for $\alpha \in A$ and \newline

$\qquad \qquad g \in {}^\kappa\text{Ord}$ we
have:
{\roster
\itemitem{ $(\alpha)$ }  for $\alpha < \beta$ in $A,t_\beta \backslash
t_\alpha \in I$ but $t_\alpha \backslash t_\beta \notin I$ \newline
(i.e. $\langle t_\alpha/I:\alpha \in A \rangle$ is strictly decreasing in 
${\Cal P}(\kappa)/I$)
\itemitem{ $(\beta)$ }  $t_\alpha = \{i < \kappa:f_\alpha(i) \le g(i)\}$.
\endroster}
\endroster
\endproclaim

\proclaim{\stag{1.6B}Claim}  1) For $\langle f_\alpha:\alpha < \delta
\rangle,\kappa,I$ as in the hypothesis of Claim \scite{1.2} such that
cf$(\delta) > \kappa^+$, the following are equivalent:
\medskip
\roster
\item "{$(a)$}"  there are $A \subseteq \delta$ unbounded, and $s_\alpha \in
I$ for $\alpha \in A$ such that: \newline
$\langle f_\alpha(i):\alpha \in A,i \in \kappa \backslash s_\alpha \rangle$
is strictly increasing in $\alpha$ for each $i < \kappa$
\item "{$(b)$}"  $(**)_I$ of Claim \scite{1.2} holds for some $f$ and
$$
\{i < \kappa:\text{cf}[f(i)] \ne \text{ cf}(\delta)\} \in I
$$
(note: $f/I$ is unique as a $<_I$-lub of $\langle f_\alpha:\alpha < \delta
\rangle$).
\endroster
\medskip

\noindent
2) If cf$(\delta) > \text{ gen}(I)$, in clause $(a)$ above without loss of
generality $s_\alpha = s$ for $\alpha \in A$; if $I = J^{bd}_\kappa$ without
loss of generality $s_\alpha = s = [i(*),\kappa)$ for some fixed $i(*)$. 
\endproclaim
\bigskip

\definition{\stag{1.6} Definition}  We define (for $\bar f$ and $J,\kappa,
\lambda,\mu$): \newline
1) $S^{gd}$ is the set of $\bar f$-good ordinals.  An ordinal
$\delta$ is weakly $\bar f$-good if $\delta$ is a limit ordinal, $\delta
< \ell g(\bar f)$ is of cofinality $> \kappa$ but $< \lambda$ such that for
some $t \subseteq \kappa$, $t \notin J$ and for some unbounded $A \subseteq
\delta$ and $\langle s_\alpha:\alpha \in A \rangle$ such that
$s_\alpha \in J$ and
for each $i < \kappa$, $\langle f_\alpha(i):\alpha \in A,i \in t \backslash
s_\alpha \rangle$ is strictly increasing (we may say $\delta$ is $\bar f$-
good by $t$); so $\bar f \restriction \delta$ has a $<_J$-eub (see below)
$f^*_\delta,f^*_\delta(i) = \sup\{f_\alpha(i):\alpha \in A \text{ and }
i \in s_\alpha\}$.  An ordinal $\delta$ is $\bar f$-good if we could
have $t = \kappa$. \newline
2) $S^{bd}$, a set of bad points for $\bar f$, is the set of limit 
ordinals $\delta < \ell g(\bar f)$
of cofinality $> \kappa$ but $< \lambda$ such that $\bar f \restriction
\delta$ satisfies $(**)_J$ of Claim \cite{Sh:355, {1.2}} - for some 
$f^*_\delta$ (i.e. $\bar f \restriction \delta$ has as a 
$<_J-eub \, f^*_\delta$ (i.e. $<_J$-exact upper bound, this means (here):
$\alpha < \delta \Rightarrow f_\alpha <_J f_\delta \text{ and } g <_J
f_\delta \Rightarrow \dsize \bigvee_{\alpha < \delta} g <_J f_\delta$)
such that $\{ i:\text{cf}[f^*_\delta(i)] \le \kappa\}
\in J)$, but $\delta \notin S^{gd}$; without loss of generality
$f^\ast_\delta \in \dsize \prod_{i<\kappa} \lambda_i$ (as $f^\ast_\delta
\le_J f_\delta$). \newline
3) $S^{ch}$, a set of $\delta$ chaotic points for $\bar f$, is the 
set of limit ordinals $\delta < \ell g(\bar f)$
of cofinality $> \kappa$ but $< \lambda$ such that:
\medskip
\roster
\item "{$(a)$}"  $\bar f \restriction \delta$ satisfies 
$(*)_{D_\delta}$ of Claim \cite{Sh:345, {1.2}} i.e. for some ultrafilter 
$D_\delta$ on $\kappa$ disjoint to $J$ and
$\langle a_i:i< \kappa \rangle$, $a_i$ a set of $\le \kappa$ ordinals,
for every $i < \delta$ for some $j < \delta$ (but $> i$)
and $g \in \dsize \prod_{i<\kappa}a_i$ we have $f_i/D < g/D < f_j/D$
\endroster
\noindent
or
\roster
\item "{$(b)$}"  $\bar f \restriction \delta$ satisfies $(***)_I$ of Claim
\cite{Sh:355, {1.2(iii)}}; (i.e. for some $A,g$ we have $A \subseteq \delta$ 
unbounded, $g \in {}^\kappa Ord$ and $\langle \{ i:f_\alpha(i) \le g(i) \} / 
J:\alpha \in A \rangle$ is strictly decreasing in ${\Cal P}(\kappa)/J$. (Hence
cf$(\delta) \le 2^\kappa$ and if $J$ is a maximal ideal, this does not 
occur).   
\endroster
\medskip

\noindent
4)  $S^{wg}(\bar f)$ is the set of weakly good $\delta < \mu$ satisfying
$\delta \notin
S^{ch}$; $S^{vbd}$ is the set $S^{bd} \backslash S^{wgd} \backslash S^{ch}$.
\enddefinition
\bigskip

\remark{\stag{1.7} Remark} 1) For the specific $J = J^{bd}_\kappa$ we could
redefine $S^{gd}$ as follows: for some unbounded $A \subseteq \delta$ and
$i(*) < \kappa$ for $\alpha \in A$, for each $i \in (i(*),\kappa)$
$\langle f_\alpha(i):\alpha \in A \rangle$ is strictly increasing.  We use:
$J$ is generated by a family of $< cf(\delta)$ subsets of $\kappa$. \newline
2) More accurately, we should have written $S^{gd}[\bar f]$ or
even $S^{gd}_J[\bar f,\lambda]$, (as $\kappa$ can be reconstructed),
similarly for $bd$ and $ch$. \newline
3)  $gd,bd,ch$ stands for good, bad, chaotic, respectively.
\endremark 
\bigskip

\demo{\stag{1.8} Fact}  1) $\langle S^{gd},S^{bd},S^{ch} \rangle$ is a 
partition of $\{ \delta < \ell g(\bar f):\kappa < cf(\delta) < \lambda \}$.
Also $\langle S^{wgd},S^{vbd},S^{ch} \rangle$ is a partition of that set.
\newline
2) $\delta \in S^{gd} \cup S^{bd}$ iff $\bar f \restriction
\delta$ satisfies $(*)_I$ of Claim \cite{Sh:355, {1.2}}, for $I = J$; i.e.
$\bar f \restriction \delta$ has a $<_J$-eub. \newline
3) $S^{gd} \subseteq S^{wg}$ and $S^{vbd} \subseteq S^{bd}$.
\enddemo
\bigskip

\demo{Proof}  If $\delta \in S^{gd}$, \underbar{then} 
$\bar f \restriction \delta$ satisfies $(*)_I$ of Claim \cite{Sh:355,
{1.2}} (for $I = J$), by \scite{1.6}.  Now Part 2) holds by the
definition of $S^{bd}$. Now Part (1) is immediate by Claim
\cite{Sh:345, {1.2}} (and Definition  \scite{1.6}).
\newline
2), 3)  Easy, too. \hfill$\square_{\scite{1.8}}$
\enddemo
\bigskip

\proclaim{\stag{1.9} Claim}  1) If $\delta \in S^{bd}$ \underbar{then} letting
$\lambda^\delta_i = \text{ cf}[f^*_\delta(i)]$ (see \scite{1.6}(2)), 
we have:
\medskip
\roster
%
\item "(i)"  $\{ i:\lambda^\delta_i \le \kappa$ or $\lambda^\delta_i >
\text{ cf}(\delta) \} \in J$ (so without loss of generality it is empty)
\item "(ii)" $\lambda^\delta_i < \lambda_i$ and
\item "(iii)"  $(\dsize \prod_{i<\kappa}\lambda^\delta_i,<_J)$
has true cofinality cf$(\delta)$ and
\item "(iv)"  for no $\lambda'$ is $\{ i:\lambda^\delta_i \ne
\lambda^\prime \} \in J$
\item "(v)"  If $\delta$ is not weakly good, then for every $\lambda'$ we
have $\{ i:\lambda^\delta_i = \lambda^\prime \} \in J$
\item "{$(vi)$}"  for every $\lambda' \ne \text{ cf}(\delta)$ we have $\{ i:
\lambda^\delta_i = \lambda'\} \in J$.
\endroster  
\medskip

\noindent
2)  If $J$ is $\theta$-complete, $\lambda = \chi^{+ \zeta}$,
$\zeta \le \theta$ a limit ordinal and $\kappa \le \chi$, \underbar{then}:
\medskip
\roster
%
\item "{$(i)$}"  $\delta \in S^{gd} \Rightarrow \{ i < \kappa:\lambda^\delta_i
\ne \text{ cf}(\delta) \} \in J$.
\item "{$(ii)$}"  $\delta \in S^{vd} \Rightarrow \{ i < \kappa:
\lambda^\delta_i \ge \chi \} \in J$.
\item "{$(iii)$}"  $\delta \in S^{wg} \Rightarrow \{ i < \kappa:
\lambda^\delta_i \ge \chi,\lambda^\delta_i \ne \text{ cf}(\delta) \} \in J$.
\item "{$(iv)$}"  cf$(\delta) \ge \chi \and \delta \in S^{wg} \backslash 
S^{gd} \Rightarrow \{i < \kappa:\lambda^\delta_i \ge \chi$ and 
$\lambda^\delta_i = \text{ cf}(\delta)\} \in J^+$.
\endroster
\medskip

\noindent
3) If $J$ is $\theta$-complete, $\zeta \le \theta$ (a limit
ordinal) and $\lambda = \kappa^{+ \zeta}$ then $S^{bd} = \emptyset$.
\endproclaim
\bigskip

\demo{Proof}  1) The first phrase holds as $f^\ast_\delta$ is $<_J$-exact
upper bound of $\bar f \restriction \delta$ (see Claim \cite{Sh:355, {1.2}}).  
For the rest see \cite{Sh:355, {1.6}}. \newline
2) Suppose $\delta \in S^{bd}$ is not weakly good and
$t =: \{ i:\lambda^\delta_i > \chi \}$ is not in $J$;
as $\delta \in S^{bd}$ trivially $\{ i:\lambda^\delta_i > \text{ cf}
(\delta) \} \in J$, so $t_1 =: \{ i:\chi < \lambda^\delta_i \le \text{ cf}
(\delta) \}$ is not in $J$.  As cf$(\delta) < \chi^{+ \zeta},\zeta \le 
\theta$, the set $\{ \sigma:\chi \le \sigma = \text{ cf}(\sigma) <
\chi^{+ \zeta}\}$ has $< |\zeta| \le \theta$ members but $J$ is 
$\theta$-complete, hence
for some $s \subseteq t_1$, $s \notin J$ and $\langle \lambda^\delta_i:i
\in s \rangle$ is constant, say $\lambda^\ast$, now also
$(\dsize \prod_{i \in s} \lambda^\delta_i,<_J)$ has true
cofinality cf$(\delta)$ (by \scite{1.9}(1)(iii)), but this implies 
$\delta$ is weakly $\bar f$-good by ``$\Leftarrow$'' of \cite{Sh:355,
{1.6(1)}}. So we have proved (ii).  For (i), (iii), (iv) the proof is similar.
(In fact for clause $(i)$ the assumption of part (2) is not needed). \newline
3) If $\delta \in S^{bd}$, by \scite{1.9}(2)(i) we have 
$\{ i < \kappa:\lambda^\delta_i > \kappa \} \in J$, so one of the 
clauses of Definition \scite{1.6}(2) fails by (ii), (iii) of part (2) of our
Claim, so 
$\delta \notin S^{bd}$, contradiction.  \hfill$\square_{\scite{1.9}}$
\enddemo
\bigskip

\proclaim{\stag{1.10} Claim}  If $\delta \in S^{gd}$, \underbar{then} 
for some club $C$ of $\delta$ we have:

$$
\alpha \in C \and \text{ cf}(\alpha) > \kappa \Rightarrow \alpha \in
S^{gd}.
$$
\medskip

\noindent
Similarly for $S^{wgd}$ (by $t$).
\endproclaim
\bigskip

\demo{Proof}  Straightforward by its definition.
\enddemo
\bigskip

\proclaim{\stag{1.11} Claim}  If $\delta \in S^{ch}$, then for 
some club $C$ of $\delta$

$$
\alpha \in C \and cf(\alpha) > \kappa \Rightarrow \alpha \in S^{ch}.
$$
\endproclaim
\bigskip

\demo{Proof}  If $\delta \in S^{ch}$, then \newline
\underbar{Case a}:  For some ultrafilter $D$ on $\kappa$ disjoint to
$J$ (say $D_\delta$) and club $C$ of $\delta$, there are $\langle g_\alpha:
\alpha \in C \rangle$ such that:
\medskip
\roster
\item "{$(a)$}"  for $\alpha < \beta$ in $C$ we have $f_\alpha/D 
< g_\alpha/D < f_\beta/D$
\item "{$(b)$}"  for each $i < \kappa$ the set $s_i = \{ g_\alpha(i):
\alpha \in C \}$ has power $\le \kappa$.
\endroster
\medskip

\noindent
Now for every $\delta^\prime \in C$ such that $\delta^\prime =
\text{ sup}(\delta^\prime \cap C)$ and $cf(\delta^\prime) > \kappa$ clearly
\newline
$\langle g_\alpha:\alpha \in C \cap \delta^\prime \rangle$, $C \cap \delta',
\langle s_i:i < \kappa \rangle$ witness that $\langle f_\alpha:
\alpha < \delta' \rangle$ also satisfies $(*)_D$ of \cite{Sh:355,
{Claim 1.2}},  hence $\delta' \in S^{ch}$. \newline
\underbar{Case b}:  (Not case (a)).  (See \scite{1.6}(3)).  Also 
trivial proof. \hfill$\square_{\scite{1.11}}$
\enddemo    
\bigskip

\proclaim{\stag{1.12} Claim}  For every regular $\chi,\kappa < \chi < 
\lambda,S^{gd}_\chi =: \{ \delta \in S^{gd}:\text{cf}(\delta) = \chi \}$ 
is stationary.
\endproclaim
\bigskip

\demo{Proof}  Suppose $C^*$ is a club of $\mu$.  We choose by induction
on $\zeta < \chi$, $\alpha_\zeta$, $g_\zeta$, $s_\zeta$ such that:
\medskip
\roster
\item  $\alpha_\zeta \in C^\ast$
\item  $\dsize \bigcup_{\xi < \zeta} \alpha_\xi < \alpha_\zeta$
\item  $g_\zeta \in \dsize \prod_{i < \kappa} \lambda_i$, $g_\zeta(i) =:
\text{ sup } \{f_{\alpha_\xi}(i) + 1:\xi < \zeta \}$
\item  $s_\zeta \in J$
\item  $i \in \kappa \backslash s_\zeta \Rightarrow g_\zeta(i) <
f_{\alpha_\zeta}(i)$.
\endroster
\medskip

In stage $\zeta$, first define $g_\zeta$ by (3) then find $\beta_\zeta <
\mu$ such that $g_\zeta <_J f_{\beta_\zeta}$ (possible as $\langle \lambda
_j:j < \kappa \rangle$ is the $<_J$-eub of $\langle f_\alpha:\alpha < \mu
\rangle$) then choose $\alpha_\zeta$, $\beta_\zeta \cup \dsize \bigcup
_{\xi < \zeta} \alpha_\xi < \alpha_\zeta \in C^*$ (possible as
$|\zeta| < \lambda$, $C^*$ unbounded in $\mu$). Lastly choose $s_\zeta$
to satisfy (5) by the definition of $<_J$.  Clearly (1)-(5) holds.

Let $\delta = \dsize \bigcup_{\zeta < \mu} \alpha_\zeta$, clearly $\delta
\in C^\ast$ (being closed) and $cf(\delta) = \chi$ and clearly $\delta \in
S^{gd}_\chi$.  [Note: if $gen(J) < \chi$, as $\chi = \text{ cf}(\chi)$ 
for some $t \in J$, $\{ \zeta < \mu:s_\zeta \subseteq t \}$ is unbounded 
in $\chi$.
Also it is natural to deduce it from the existence of a stationary
\newline
$S \subseteq \{ \delta < \mu:\text{cf}(\delta) = \chi\},S \in I[\mu]$ (see
\cite{Sh:420, {\S1}})].
\enddemo
\bigskip

\noindent
We shall not use
\proclaim{\stag{1.13} Lemma}  Suppose $\lambda = \kappa^{+ \kappa},\langle 
\lambda_i:i < \kappa \rangle,\mu,J$ are as in $(*)$ of \scite{1.4} and in 
addition $S^{bd} = \emptyset$.  Let $C^* = \{ \delta < \lambda^+:\delta
\text{ divisible by } \lambda^2 \}$.  Then we can find $\langle \alpha^\delta
_{i,j}:\delta \in C^* \cap S^{gd},i < \text{ cf}(\delta),j < \kappa \rangle$ 
such that:
\medskip
\roster
\item "{$(a)$}"  $\alpha^\delta_{i,j} < \delta$ (and, if you want, $\langle
\alpha^\delta_{i,j}:i < \text{ cf}(\delta) \rangle$ is increasing continuous 
with limit $\delta,\langle \alpha^\delta_{i,j}:j < \kappa \rangle$ increasing 
with $j$).
\item "{$(b)$}"  for every $\delta(*) < \lambda^+$ we can find 
$\langle C^\delta,s^\delta_i:\delta \in \delta(*) \cap C^* \cap S^{gd},
i < cf(\delta) \rangle$ such that:
{\roster
\itemitem{ ($\alpha$) }  $C^\delta$ a club of cf$(\delta)$
\itemitem{ ($\beta$) }   $s^\delta_i \in J$
\itemitem{ ($\gamma$) }   $\langle \{ \alpha^\delta_{i,j}:i \in C^\delta,
j \notin s^\delta_i \}: \delta \in \delta (*) \cap C^* \cap S^{gd} \rangle$ 
are pairwise disjoint.
\endroster}
\endroster
\endproclaim

\noindent
This we refer to \cite{MgSh:204}, anyhow see proof of \scite{2.5}. 

\newpage

\head {\S2 Nicely Compact Set and Non-Reflection} \endhead
\resetall
\bigskip

This section answers the following two questions.  \newline
One question of Foreman and Magidor asks for consistency of the
following reflection principle:
\medskip
\roster
\item "${(*)_\lambda}$"  \underbar{if} $S \subseteq {\Cal S}_{<\aleph_2}
(\lambda)$ is a stationary set (of ${\Cal S}_{<\aleph_2}(\lambda)$, not of
$\lambda$), each $a \in S$ is $\omega$-closed (i.e. cf$(\delta) = \aleph_0
\and \delta = \text{ sup }(\delta \cap a) < \lambda \Rightarrow \delta \in
a)$ \underbar{then} for some $B \subseteq \lambda$, $|B| = \aleph_2$ and
$S \cap {\Cal S}_{<\aleph_2}(B)$ is stationary.
\endroster
\medskip

\noindent
We shall show that it is very hard to get them (in a way defeating the
application they have in mind: $\theta < \aleph_2$ follows from large
cardinals). E.g. if $\mu$ is strong limit singular of cofinality $\aleph_0$,
$2^\mu = \mu^+$ then $(*)_{\mu^+}$ is false.

Another question of 
Mekler and Eklof, improve somewhat \cite{MgSh:204} showing in particular 
that if we have gotten in the (regular) $\lambda,\kappa$ examples for
incompactness then we can get one for $\lambda^{+(\kappa+1)}$.

To do this (see below) if $\lambda > \kappa^{+ \kappa}$, we need a
counterexample for $\lambda$, which without loss of generality is as in
\cite{Sh:161} canonical examples in which only cardinals in which examples 
exist appear (see below).  For this we
define when a set $K$ of cardinals is nicely incompact and prove that the
family of such sets is closed under the required operation.
\bigskip   

\definition{\stag{2.1} Definition}  A set $K$ of 
regular cardinals is called nicely incompact if:
\medskip
\roster
\item "{$(a)$}" for each $\lambda \in K$ there is a $\lambda$-free not
$\lambda^+$-free abelian group or $\lambda = \aleph_0$
\item "{$(b)$}"  \underbar{if} $\lambda$ is in $K,\kappa_1,\cdots,\kappa_m
\in K \cap \lambda$ \underbar{then} there is an example \newline
$\langle s^\ell_\eta:\eta \in {\frak S}_f,\ell < m \rangle$ of incompactness
(see \cite{Sh:161, {\S3}} or \cite{Sh:521}, we assume knowledge of e.g.
\cite{Sh:521, {AP}}) such that:

$$
\align
&(\lambda_{<>} = \lambda \text{ and) for some } \langle k_1,\cdots,k_m
\rangle \text{ we have} \\
& \eta \in S \and \ell g(\eta) = k_\ell \Rightarrow \lambda_\eta = 
\kappa_\ell.
\endalign
$$
\endroster
\enddefinition 
\bigskip

\proclaim{\stag{2.2} Claim}  1) If $K$ is nicely incompact, 
$\lambda = \text{ max}(K),\lambda$ regular \underbar{then} $K \cup 
\{ \lambda^+ \}$ is nicely incompact. \newline
2)  Assume $K = \{ \lambda_i:i < \delta \}$ increasing and $j < \delta
\Rightarrow \{ \lambda_i:i \le j \}$ nicely incompact.  \underbar{Then} 
$K$ is nicely incompact. \newline
3)  $\{ \aleph_0 \}$ is nicely incompact.
\endproclaim
\bigskip

\demo{Proof}  Easy.
\enddemo
\bigskip

\remark{Remark}  If proof of \scite{2.2}(1) is not clear - read proof of 
\scite{2.3} and throw away most (use $\{ \delta < \lambda^+:\text{cf}(\delta) 
= \lambda \}$ here
instead $\{ \delta \in A_g:\text{cf}(\delta) = \lambda'\}$ there.
\endremark
\bigskip

\proclaim{\stag{2.3} Theorem}  If $K$ is nicely incompact, $\chi,\kappa \in 
K,K$ has no last element and $\chi^{+ \kappa} \ge \text{ sup } (K)$ 
\underbar{then} $K \cup \{ \chi^{+(\kappa +1)} \}$ is nicely incompact.
\endproclaim
\bigskip

\remark{Remark}  By the closure properties \scite{2.2}(1), \scite{2.2}(2), 
\scite{2.3} (and as a starting point \scite{2.2}(3): $\{ \aleph_0 \}$ 
being nicely incompact) we get a family of $\aleph_0$ cardinals which 
probably is the minimal set of incompactness.
\endremark
\bigskip

\demo{Proof}  If $\chi < \kappa = \aleph_\kappa$, then $\chi^{+(\kappa +1)}
= \kappa^+$ and the desired conclusion holds is by Claim \scite{2.2}(1).

If $\chi < \kappa < \aleph_\kappa$ then $\chi^{+\kappa} = \left( \kappa^+
\right) ^\kappa = \aleph_\kappa$ so without loss of generality
$\chi = \kappa^+$.

So without loss of generality $\chi \ge \kappa^+,\lambda = \chi^{+\kappa},
\mu = \lambda^+,\lambda > \kappa$ and $J$, \newline
$\langle \lambda_i:i < \kappa 
\rangle,\bar f = \langle f_\alpha:\alpha < \lambda \rangle$ be as in $(*)$ 
of \scite{1.4} and for simplicity \newline
$f_{\alpha_1}(i_1) = f_{\alpha_2}(i_2) \Rightarrow i_1 = i_2$, e.g.
$f_\alpha(i) > \dsize \sum_{j < i} \lambda_j$; and let $S^{gd},S^{bd},S^{ch}$
be from Definition \scite{1.6}.  Remember that by \scite{1.10},
\scite{1.11} we have
\medskip
\roster
\item "$(*)$"  if $x \in \{ gd, wgd,ch\},\delta \in S^x$ \underbar{then} 
for some club $E$ of $\delta$ we have: \newline
$\alpha \in E \and \text{ cf}(\alpha) \ge \kappa^+ \Rightarrow \alpha \in
S^x$.
\endroster
\medskip

\noindent
In the definition of ``$K$ is nicely incompact" we have to check only for
$\mu = \lambda^+$.
Let $m < \omega$ and $\kappa_1,\dotsc,\kappa_m \in K \cap \lambda^+$.  So 
$\langle \kappa_1,\dotsc,\kappa_m \rangle \in K \cap \chi^{+\kappa}$,
without loss of generality $\kappa_1 < \kappa_2 < \dotsc < \kappa_m$. \newline
Let $\lambda' \in K$ be such that $\kappa < \lambda,\chi < \lambda'$ and
$\kappa_m < \lambda'$ (there is such a $\lambda'$ as $K$ has no last member).

Let $S = \{ \delta \in S^{gd}(\bar f):\text{cf}(\delta) = \lambda'\}$ (by 
\scite{1.12} it is stationary) and for $\delta \in S$ choose 
$E_\delta \subseteq \delta$ an unbounded subset of order type $\lambda'$.

Let $E_\delta = \{ \alpha^\delta_\zeta:\zeta < \lambda^\prime \}$
(increasing). \newline
Now choose an example for $\lambda^\prime:\langle s^\ell_\eta:\eta \in 
{\frak S}'_f \text{ and } \ell < n(*)\rangle,{\frak S}'$ a 
$\lambda'$-set (see \cite{Sh:161, {\S3}} or \cite{Sh:521, {AP}}), 
such that for some sequence
$\langle k_1,\dotsc,k_n,k \rangle$ of natural number $\le lg(\eta)$ for
$\eta \in {\frak S}'$

$$
\gather
\eta \in {\frak S}^\prime \and \ell g(\eta) = k_i
\Rightarrow \lambda_\eta = \kappa_\ell \\
\eta \in {\frak S}^\prime \and \ell g(\eta) = k \Rightarrow
\lambda_\eta = \kappa.
\endgather
$$
\medskip

\noindent
(possible by Definition \scite{2.1}). \newline
Lastly, let

$$
{\frak S} = \{ <> \} \cup \left\{ <\delta> \char 94 \eta:
\delta \in S \subseteq \lambda^+ \text{ and }
\eta \in {\frak S}' \right\}.
$$
\medskip

\noindent
For $\eta \in {\frak S}_f$ let (denoting $\delta = \eta(0)$ which
necessarily is in $S$ and letting \newline
$\eta = \langle \delta \rangle \char 94
\nu$ so $\nu = \eta \restriction [1,\ell g(\eta))$):

$$
s^{\ell + 1}_\eta = (s^\ell_\nu \times \{ \delta \}) \text{ if }
\ell < n(*), \text{ and } s^0_\eta = \left\{ \langle f_\delta
(\eta(k)),\eta \restriction [1,\ell g(\eta)),\alpha^\delta_{\eta(k)})\rangle:
k < \omega \right\}.
$$

$$
A_\eta = \cup \{ s^\ell_\eta:\ell \le n(*) \}
$$
\medskip

First note that $\langle A_\eta:\eta \in {\frak S}_f \rangle$ has the
right form hence has no one-to-one choice function (see \cite{Sh:161,
{\S3}} or \cite{Sh:521, {AP}}).

Secondly, to prove every subfamily of cardinality $< \mu$ has a one-to-one
choice function it is enough to prove \scite{2.5} below.  It applies with
$\lambda',\mu$ here standing for $\chi,\mu^*$ there.  The no reflection is
by \scite{2.4} below.
\enddemo
\bigskip

\demo{\stag{2.4} Observation}  Let $\bar \lambda,J,\mu,\lambda$ be as in
$(*)$ of \scite{1.4}. \newline
1) If $J$ is $\theta$-complete, $\lambda = \chi^{+ \zeta},\zeta \le \theta$
a limit ordinal, cf$(\chi) > \kappa$, then:
\medskip
\roster
%
\item "{$(i)$}"  $\delta \in S^{gd} \Rightarrow \{ i < \kappa:
\lambda^\delta_i \ne \text{ cf}(\delta)\} \in J$
\item "{$(ii)$}"  $\delta \in S^{vbd} \Rightarrow \{ i < \kappa:
\lambda^\delta_i \ne \text{ cf}(\delta)\} \in J$
\item "{$(iii)$}"  $\delta \in S^{wg} \Rightarrow \{ i:\lambda^\delta_i \ge 
\chi,\lambda^\delta_i \ne \text{ cf}(\delta)\} \in J$
\item "{$(iv)$}"  cf$(\delta) \ge \chi \and \delta \in S^{wg} \backslash
S^{gd} \Rightarrow \{ i < \kappa:\lambda^\delta_i \ge \chi$ and
$\lambda^\delta_i = \text{ cf}(\delta)\} \in J^+$.
\endroster
\medskip

\noindent
2) If $J$ is $\theta$-complete, $\zeta \le \theta$ (a limit ordinal) and
$\lambda = \chi^{+ \zeta}$ then for no $\delta \in S^{bd}$ do we have:
$\{ \alpha \in S^{gd}:\alpha < \delta,\text{cf}(\alpha) \ge \chi\}$ is
stationary.
\enddemo
\bigskip

\demo{Proof}  1) Like \scite{1.9}. \newline
2)  Let $\delta$ be a counterexample, let $f^*_\delta$ be a
$<_J$-eub of $\bar f \restriction \delta$ and by part (1) without loss of
generality $s^* = \{i < \kappa:\lambda^i_\delta = \text{ cf}[f^*_\delta(i)]
< \chi\} \} \in J^+$.  Choose $E_i \subseteq f^*_\delta(i)$ a club of
order type $\lambda^i_\delta$. \newline
Choose $\langle \alpha_\zeta:\zeta < \text{ cf}(\delta) \rangle$ is a
strictly increasing sequence of ordinals with limit $\delta$.  Now we can
choose by induction on $\zeta < \text{ cf}(\delta)$ a pair $(\beta_\zeta,
g_\zeta)$ such that:
\medskip
\roster
\item "{$(a)$}"  $g_\zeta \in \dsize \prod_{i < \kappa} E_i$ such that
\newline
$\varepsilon < \zeta \Rightarrow f_{\beta_\varepsilon} \restriction s^* <_J
g_\zeta \restriction s^*$
\item "{$(b)$}"  $\beta_\zeta \in \left( \dsize \bigcup_{\varepsilon <\zeta}
\beta_\varepsilon \cup \alpha_\zeta,\delta \right)$ is such that
$g_\zeta \restriction s^* <_J f_{\beta_\zeta} \restriction s^*$.
\endroster
\medskip

\noindent
So if $\delta$ is a counterexample, then for some $\zeta < \text{ cf}(\delta)$
we have: $\text{cf}(\zeta) \ge \chi$ and $\beta^* = \dsize 
\bigcup_{\varepsilon < \zeta} \beta_\varepsilon \in S^{gd}$.  Hence there 
is an increasing sequence $\langle \gamma_j:j < \text{ cf}(\beta^*) \rangle$ 
of ordinals with limit $\beta^*$,
and $t_j \in J$ for $j < \text{ cf}(\beta^*)$ such that for each $i,
\langle f_{\gamma_j}(i):j < \text{ cf}(\beta^*),i \notin t_j \rangle$ is
strictly increasing.  For each $i \in s^*$ there is $j(i) < \text{ cf}
(\beta^*)$ such that all the ordinals $\{f_{\gamma_j}(i):j(i) \le j <
\text{ cf}(\beta^*) \text{ and } i \notin t_j\}$ realize the same
Dedekind cut of $E_i$.  Let $j(*) = \underset{i < \kappa} \to \sup j(i) <
\text{ cf}(\beta^*)$ (as cf$(\beta^*) \ge \text{ cf}(\chi) > \kappa$).  Now
the contradiction should be clear. \hfill$\square_{\scite{2.4}}$
\enddemo
\bigskip

\remark{\stag{2.4A} Remark}  1) Alternatively, in the proof of \scite{2.3} 
redefine $S^{ch}(\bar f)$ replacing \newline
``$|s_i| \le \kappa$" by ``$|s_i| \le \chi$".  We call
it $S^{ch}_\chi(\bar f)$ (similarly the others). \newline
2)  Note: if $\lambda$ belongs to some nicely incompact $K^*$, then there is 
a finite nicely incompact $K \subseteq K^*$ to which $\lambda$ belongs.
\endremark
\bigskip

\proclaim{\stag{2.5} Claim}  In the context of $(*)$ of \scite{1.4}, suppose
$\kappa < \chi = \text{ cf}(\chi) < \mu^* \le \mu$ and 
$S \subseteq \{ \delta \in S^{gd}:\text{cf}(\delta) = \chi \}$ is 
stationary but reflect in no $\delta \in S^{bd}$ of cofinality $< \mu^*$ and
$f_{\alpha_1}(i_1) = f_{\alpha_2}(i_2) \Rightarrow i_1 = i_2$ (e.g.
$f_\alpha(i) > \underset{j < i} \to \sup \lambda_j$).  
For $\delta \in S$ let $E_\delta$ be
an unbounded subset of $\delta$.  Let $A_\delta = E_\delta \times
(\text{Rang }f_\delta) \subseteq \delta \times \lambda$.  \underbar{Then}
$\langle A_\delta:\delta \in S \rangle$ has no one-to-one choice function
but is $\mu^*$-free in the following sense:
\medskip
\roster
\item "$(*)$"  if $B \subseteq S$, $|B| < \mu^\ast$ then we can find
$\alpha_\delta < \delta$ and $s_\delta \in J$ for $\delta \in B$ such that
the sets $\langle (E_\delta \backslash \alpha_\delta) \times \text{ Rang }
(f_\delta \restriction (\kappa \backslash s_\delta):\delta \in B \rangle$
are pairwise disjoint.
\endroster
\endproclaim
\bigskip

\remark{Remark}  1) When we use it in \scite{2.3}, $\mu^* = \mu$ and 
the non-reflection is guaranteed by \scite{2.4}. \newline
2) If we want to replace $A_\delta$ by a subset of $\delta$,
assume every $\delta \in S$ is divisible by $\lambda^2$, and every $\alpha
\in E_\delta$ is divisible by $\lambda$, and let

$$
A_\delta = \{ \alpha^\delta_\zeta + f_\delta(i):i < \kappa,
\theta < \chi \} \text{ where }
S_\delta = \{ \alpha^\delta_\zeta:\zeta < \chi \}
\text{ (increasing) }.
$$
\medskip

\noindent
3) If $J = J^{bd}_\kappa$ it may look nicer to assume
$\langle \text{Rang } f_\alpha \restriction[i,\kappa):\alpha \in E_\delta
\rangle$ (for suitable $i = i_\delta < \kappa$) are pairwise disjoint (not
hard to get as $\delta \in S^{gd}$).  Then use $f_\alpha \restriction [i,
\kappa)$ only. \newline
4)  It seems reasonable to demand (and it causes no problem)
$i < \kappa \and \alpha < \mu \Rightarrow f_\alpha(i) > {\text{sup}} \{
f_\beta(j):j < i,\beta < \mu \}$.
\endremark
\bigskip

\demo{Proof}  This is proved by induction on $\cup \{ \delta + 1:\delta \in
B \}$. \newline
\noindent  Without loss of generality $B \ne \emptyset$.
\enddemo
\bigskip

\demo{Case 1}  $B$ has a last element $\delta(*)$.

By the induction hypothesis there are $\langle (\alpha^1_\delta,
s^1_\delta):\delta \in B \backslash \{ \delta(*) \} \rangle$ as required for
$B \backslash \{ \delta(*) \}$.  For each $\delta \in B \cap \delta(*)$ we
have $t_\delta = \{ i < \kappa:f_\delta(i) \ge f_{\delta(*)}(i) \} \in J$.
Let us define for $\delta \in B$:
$$
\alpha_\delta =
\left\{
\aligned \alpha^1_\delta \quad &\text{\underbar{if}} \quad \delta \in B,
\delta \ne \delta(*)\\
0 \quad &\text{\underbar{if}} \quad \delta = \delta(*)
\endaligned
\right.
$$
\medskip
$$
s_\delta =
\left\{
\aligned s_\delta \cup t_\delta \quad &\text{\underbar{if}} \quad
\delta \in B, \delta \ne \delta(*)\\
\emptyset \quad &\text{\underbar{if}} \quad \delta = \delta(*).
\endaligned
\right.
$$
\enddemo
\bigskip

\demo{Case 2}  $B$ has no last element, $\delta(*) = \text{ sup}(B)$ has
cofinality $\le \chi$.

Let $\langle \gamma_\zeta:\zeta < cf\delta(*) \rangle$ be increasing
continuous sequence with limit $\delta(*),\gamma_0 = 0$ and
cf$(\gamma_\zeta) \ne \chi$, hence $\gamma_\zeta \notin S$.  For each
$\zeta$, apply the induction hypothesis to $B_\zeta =: B \cap (\gamma_\zeta,
\gamma_{\zeta +1})$ and get $\langle (\alpha^\zeta_\delta,s_\delta):\delta
\in B_\zeta \rangle$.

We let for $\delta \in B$ (note: $B$ is the disjoint union of $B_\zeta$
for $\zeta < \chi$):
$$
\alpha_\delta = {\text{Max}} \{ \alpha^\zeta_\delta,\gamma_\zeta \}
{\text{ for }} \delta \in B_\zeta
$$
$$
s_\delta = s^\zeta_\delta {\text{ for }} \delta \in B_\zeta.
$$
\medskip

\noindent
Easy to check that they are as required.
\enddemo
\bigskip

\demo{Case 3}  $B$ has no last element, $\delta(*) = \text{ sup}(B)$ has
cofinality $> \chi$ and $S^{gd}$ (hence $S$) does not reflect in 
$\delta(*)$. \newline
Same proof as Case 2.
\enddemo
\bigskip

\demo{Case 4}  $B$ has no last element, $\delta(*) = \text{ sup}(B)$ has
cofinality $> \chi$ and $S^{gd}$ reflects in $\delta(*)$.

As $|B| < \mu^\ast$, clearly $cf(\delta(*)) < \mu^\ast$ hence (by assumption)
$\delta(*) \notin S^{bd}$.  Also $\delta(*) \notin S^{ch}$ (by \scite{1.11} 
as: $S^{gd} \cap \delta(*)$ is stationary and $S^{gd} \cap S^{ch} = \emptyset$
(which holds by \scite{1.8})).  So $\delta(*) \in S^{gd}$, hence we can find
$\langle \gamma_{\zeta +1}:\zeta < \text{ cf}(\delta(*)) \rangle$ 
increasing with limit $\delta(*)$, and $s_\zeta \in J$ for $\zeta < 
\text{ cf}(\delta(*))$ be such that:
\medskip
\roster
\item "{$(*)$}"  $i < \kappa \and \zeta < \xi < \text{ cf}(\delta(*)) 
\and i \notin s_\zeta \cup
s_\xi \Rightarrow f_{\gamma_\zeta}(i) < f_{\gamma_\xi}(i)$
\item "{$(**)$}"  $\gamma_\zeta \in S \Rightarrow \text{ cf}(\zeta) = \chi$.
\endroster
\medskip

\noindent
For limit $\zeta < \text{ cf}(\delta(*))$ let $\gamma_\zeta =: \dsize \bigcup
_{\xi < \zeta}\gamma_{\xi + 1}$ and let $\gamma_0 = 0$.

Now if $\zeta < \text{ cf}(\delta(*))$ is limit and $\gamma_\zeta \in S$ (so
cf$(\zeta) = \chi$), let

$$
t^0_\zeta =: \left\{
i:\text{ for arbitrarily large } \xi < \zeta, i \in \kappa \backslash
s_{\xi + 1} \text{ and } f_{\gamma_{\xi + 1}}(i) \ge f_{\gamma_\zeta}
(i) \right\}
$$
\medskip

\noindent
now $t^0_\zeta$ is in $J$ [otherwise, for each $i \in t^0_\zeta$ 
let $Y^i_\zeta = \{ \xi < \zeta:i \in \kappa \backslash s_{\xi +1}$ and
$f_{\gamma_{\xi + 1}}(i) \ge f_{\gamma_\zeta}(i)\}$ by the choice of
$\langle s_{\xi +1}:\xi < \text{ cf}(\delta(*))\rangle$, clearly

$$
\xi_1 < \xi_2 < \zeta \and \xi_1 \in Y^i_\zeta \and i \in \kappa \backslash
s_{\xi_2} \Rightarrow \xi_2 \in Y^i_\zeta.
$$
\medskip

\noindent
Let $\xi(i) \le \zeta$ be Min$(Y^i_\zeta \cup \{ \zeta\})$, and
$\xi(*) = \text{ Min}\{\xi(i):i < \kappa,\xi(i) < \zeta\}$, now $\xi(*) <
\zeta$ as cf$(\zeta) = \chi > \kappa$; and look at what occurs for
$\xi = \xi(*) + 1$].

Let $t^1_\xi = \{ i:\text{ for every large enough } \xi < \zeta,i \notin 
s_{\xi + 1} \}$.  Again easily $t^1_\xi \in J$.  Lastly let
$t^2_\xi = \{i:f_{\gamma_\xi}(i) \ge f_{\gamma_{\xi + 1}}(i)\}$ now as
$\gamma_\xi < \gamma_{\xi + 1}$ clearly $t^2_\xi \in J$ and $t^3_\zeta =
s_{\zeta + 1} \in J$.

Now for each $\zeta < \text{ cf}(\delta(*)),B_\zeta =: B \cap (\gamma_\zeta,
\gamma_{\zeta+1})$ satisfies the induction hypothesis, hence we have
$\langle (\alpha^\zeta_\delta,s^\zeta_\delta):\delta \in B_\zeta \rangle$
as required.

Notice that $B$ is partitioned to $B_\zeta$ (for $\zeta < \text{ cf}
(\delta(*))$ and \newline
\noindent $\{ \gamma_\zeta:\zeta \text{ limit } < \text{ cf}(\delta(*)) \} 
\cap S$.
We define
\medskip

$$
\alpha_\delta =
\left\{
\aligned
{\text{Max}} \{ \alpha^\zeta_\delta,\gamma_\zeta \} \quad
 &\text{\underbar{if}} \quad \delta \in B_\zeta,\, \zeta < cf(\delta(*))\\
0 \quad &\text{\underbar{if}} \quad \delta \in
\{ \gamma_\zeta:\zeta {\text{ limit }} \} \cap S 
\endaligned
\right.
$$
\medskip

\noindent
for $\zeta < \text{ cf}(\delta(*))$ and $\delta \in B_\zeta$ let

$$
t_{\zeta,\delta} = s_\zeta \cup s_{\zeta +1} \cup
\{i < \kappa:\neg[f_{\gamma_\zeta}(i) < f_\delta(i) < f_{\gamma_{\zeta +1}}
(i)]\} \in J
$$
\medskip

$$
s_\delta =
\left\{
\aligned
s^\zeta_\delta \cup t_{\delta,\zeta} \quad
 &\text{\underbar{if}} \quad \delta \in B_\zeta,\,
\zeta < \text{ cf}(\delta(*))\\
t^0_\zeta \cup t^1_\zeta \cup t^2_\zeta \cup t^3_\zeta \quad
 &\text{\underbar{if}} \quad \delta = \gamma_\zeta, \,
\zeta < \text{ cf}(\delta(*)) {\text{ limit }},\,\delta \in B(\subseteq S)
\endaligned
\right.
$$
\medskip

\noindent
It is easy to check that $\langle (\alpha_\delta,s_\delta):\delta \in B
\rangle$ is as required: let $\delta_1,\delta_2 \in B$ be distinct, without
loss of generality $\delta_1 < \delta_2$, now 
\medskip

\noindent
\underbar{Case 1}:  If $\delta_1,\delta_2 \in B_\zeta$ for some $\zeta$;
\underbar{then} use the choice of $\langle \alpha^\zeta_{\delta_1},
s^\zeta_{\delta_1} \rangle,\langle \alpha^\zeta_{\delta_2},s^\zeta_{\delta_2}
\rangle$.
\bigskip

\noindent
\underbar{Case 2}:  If $\delta_1 \in B_{\zeta_1},\delta_2 \in B_{\zeta_2},
\zeta_1 \ne \zeta_2$ \underbar{then} note $\alpha_{\delta_2} \ge
\gamma_{\zeta_2}$.
\bigskip

\noindent
\underbar{Case 3}:  If $\delta_1,\delta_2 \in \{\gamma_\zeta:\zeta <
\text{ cf}(\delta(*))\} \cap S$, \underbar{then} let $\delta_1 =
\gamma_{\zeta_1},\delta_2 = \gamma_{\zeta_2}$ hence by $(**)$ we have
$\zeta_1 + 1 < \zeta_2$, so $i \in \kappa \backslash t^2_{\zeta_1}
\Rightarrow f_{\gamma_{\zeta_1}}(i) < f_{\gamma_{\zeta_1}+1}(i)$ and
$i \in \kappa \backslash t^3_\zeta \backslash t^0_\zeta \backslash t^1_\zeta
\Rightarrow f_{\gamma_{\zeta_1}+1}(i) < f_{\gamma_{\zeta_2}}(i)$ so the
conclusion is easy.
\bigskip

\noindent
\underbar{Case 4}:  If $\delta_2 \in B_{\zeta_2},\delta_1 \in \{\gamma_\zeta:
\zeta < \text{ cf}(\delta(*))\} \cap S$ \underbar{then} use 
$\alpha_{\delta_2} \ge \gamma_{\zeta_2}$.
\bigskip

\noindent
\underbar{Case 5}:  If $\delta_1 \in B_{\zeta_1},\delta_2 =
\gamma_{\zeta_2} \in \{\gamma_\zeta:\zeta < \text{ cf}(\delta(*))\} \cap S$ 
\underbar{then} $\zeta_2 \ge \zeta_1$ hence by $(**)$ we have
$\zeta_2 > \zeta_1 + 1$, hence $i \in \kappa \backslash s_{\zeta_1 +1}
\backslash t^0_{\zeta_2} \backslash t^1_{\zeta_2} \Rightarrow
f_{\gamma_{\zeta_1}+1}(i) < f_{\gamma_{\zeta_2}}(i)$ hence
$i \in \kappa \backslash s_{\delta_1} \backslash s_{\delta_2} \Rightarrow
i \in \kappa \backslash t_{\zeta,\delta} \backslash s_{\zeta_1+1} \backslash
t^0_{\zeta_2} \backslash t^1_{\zeta_2} \Rightarrow f_{\delta_1}(i) <
f_{\gamma_{\zeta_1}+1}(i) < f_{\delta_2}(i)$.  \hfill$\square_{\scite{2.5}}$
\enddemo
\bigskip

\demo{\stag{2.6} Conclusion}  Suppose cf$(\lambda) = \kappa < \lambda,\kappa <
\chi = \text{ cf}(\chi) < \lambda$.  Then there is a stationary
\noindent  $S \subseteq \{ \delta < \lambda^+:\text{cf}(\delta) = \chi\}$ 
and $\langle A_\delta:\delta \in S \rangle$ such that:
\medskip
\roster
\item "$(\alpha)$"  $A_\delta = \left\{ \alpha^\delta_{i,\zeta}:i < \chi,
\zeta < \kappa \right\} \subseteq \delta$
\item "$(\beta)$"  $\alpha^{\delta^1}_{i_1,\zeta_1} = 
 \alpha^{\delta^2}_{i_2,\zeta_2} \Rightarrow i_1 = i_2 \and
\zeta_1 = \zeta_2$
\item "$(\gamma)$" for each $\delta \in S$ , $\zeta < \kappa$ we have
$\dsize \bigwedge_\gamma$ [the sequence $\langle \alpha_{i,\zeta}:i < \chi
\rangle$ \newline
is increasing with limit $\delta$]
\item "$(\delta)$"  $\langle A_\delta:\delta \in S \rangle$ is
$\chi^{+\kappa+1}$-free in the strong sense of \scite{2.5}; i.e. \newline
if $E \subseteq S,|E| \le \chi^{+\kappa}$ then we can find 
$\langle \langle i_\delta,\zeta_\delta \rangle:\delta \in E \rangle$ 
such that:
{\roster
\itemitem{ (a) }  $i_\delta < \chi$, $\zeta_\delta < \kappa$
\itemitem{ (b) }  the sets $A^\prime_\delta = \{ \alpha^\delta_{i,\zeta}:
\zeta_\delta < \zeta < \kappa \text{ and } i_\delta < i < \chi \}$ (for
$\delta \in E$) are pairwise disjoint.
\endroster}
\endroster
\enddemo
\bigskip

\demo{Proof}  Combine the previous claims. \newline
(Note: $\lambda$ necessarily is $\ge \chi^{+\kappa+1}$ but can be strictly
bigger).
\enddemo
\bigskip

\demo{\stag{2.7} Conclusion}  1) If $\lambda > \text{ cf}(\lambda) = \kappa,
\kappa < \chi = \text{ cf}(\chi) < \lambda$, \underbar{then} we can 
find a stationary $S \subseteq \{ \delta < \lambda^+:\text{cf}(\delta) = 
\chi\}$ and $A_\delta \subseteq \delta = \text{ sup}(A_\delta),|A_\delta|
= \chi$ for $\delta \in S$, such that:
\medskip
\roster
\item "$(*)$"  if $B \subseteq S$, $|B| \le \chi^{+\kappa}$ then
$\{ A_\delta:
\delta \in B \}$ has a one-to-one choice function.
\endroster
\medskip

\noindent
2) If in (1) $\diamondsuit_S$ holds, without loss of generality
$\{ A_\delta:\delta \in S \}$ is a stationary subset of
${\Cal S}_{<\chi}(\lambda^+)$ each $A_\delta$ is $\omega$-closed
(and even $(< \chi)$-closed).\newline
3)  $\diamondsuit_S$ holds if $2^\lambda = \lambda^+ \and \lambda \ge
\beth_\omega$ (also if $\lambda < \beth_\omega$ this is normally O.K.).
\newline
4)  Condition $(*)$ of (1) implies that for no $A^\ast \subseteq \lambda$
of power $\le \chi^{+\kappa}$, is \newline
\noindent  $\{ A_\delta:A_\delta \subseteq A^* \}$
a stationary subset of ${\Cal S}_{<\chi^+}(A^*)$.
\enddemo
\bigskip

\demo{Proof}  1) By \scite{2.6}. \newline
2)  We can replace $A_\delta$ by any $A^\prime_\delta$, $A_\delta \subseteq
A^\prime_\delta \subseteq \delta$, $|A^\prime_\delta| = \chi$, and preserve
the conclusion of (1); by $\diamondsuit_S$ we can do it as to get
stationary set. \newline
3)  By \cite{Sh:460} (see more there). \newline
4)  Check. \hfill$\square_{2.7}$
\enddemo
\bigskip

\demo{\stag{2.8} Observation}  Suppose cf$(\lambda) = \kappa < \lambda$ and 
there is $\langle A_\alpha:\alpha < \lambda^+ \rangle$, \newline
$A_\alpha \in [\lambda]^\kappa$, such that:

$$
\gather
\text{ for each } \beta < \lambda^+ \text{ we can find }
\langle A^\prime_\alpha:\alpha < \beta \rangle \text{ pairwise disjoint
with }\\
A^\prime_\alpha \subseteq A_\alpha \and |A_\alpha \backslash
A_\alpha| < \kappa.
\endgather
$$
\medskip

\noindent
Then:
\medskip
\roster
\item "{$(A)$}"  for every regular $\chi \in (\kappa,\lambda)$ we can find
\newline
$\langle \langle \alpha^\delta_{i,\zeta}:i < \chi,\zeta < \kappa \rangle:
\delta < \lambda^+,\text{cf}(\delta) = \chi \rangle$ as in \scite{2.6}, 
$\lambda^+$-free, in \scite{2.5}'s sense; i.e. see clause $(\delta)$ of
\scite{2.6}.
\item "{$(B)$}"  If $K$ is a nicely incompact set, $\kappa = \text{ max}(K),
\lambda \ge \text{ sup}(K)$ and \underbar{then} \newline
$K \cup \{ \lambda^+ \}$ is nicely incompact.
\endroster
\enddemo
\bigskip

\demo{Proof}  Trivial.  For (A) just replace every ordinal of $< \lambda$ by
$\chi$ ordinals; and copy them on each interval $[\lambda \alpha,\lambda
\alpha + \lambda)$.  Also (B) is immediate. \hfill$\square_{\scite{2.8}}$
\enddemo
\bigskip

\demo{\stag{2.9} Conclusion} If $\kappa = \text{ cf}(\lambda) < \lambda$ and
$pp_{\Gamma(\kappa)}(\lambda) > \lambda^+$ \underbar{then} the assumption 
of \scite{2.8} hence the conclusions of \scite{2.8} and \scite{2.6} holds.
\enddemo
\bigskip

\demo{Proof}  By \cite{Sh:355, {1.6(2)}} the assumptions of \scite{2.8} holds.
\enddemo
\newpage

\head {\S3 NPT is Not Necessarily Transitive} \endhead
\resetall
\bigskip

\proclaim{\stag{3.1} Lemma}  Assume the consistency of $ZFC +$ ``there are two
weakly compact cardinals $> \aleph_0$".  \underbar{Then} it is consistent 
with $ZFC$ that for some regular $\lambda > \mu > \aleph_0,NPT(\lambda,\mu) 
+ NPT(\mu,\aleph_0)$ does not imply
$NPT(\lambda,\aleph_0)$; (i.e. in some generic extension of universe $V$
this holds provided that $V \models ``\aleph_0 < \kappa < \lambda$,
$\kappa$ and $\lambda$ are weakly compact").
\endproclaim
\bigskip

\remark{\stag{3.1A} Remark}  We give more 
specific information in the three cases below.
\endremark
\bigskip

\demo{Proof}  First assume $\kappa < \mu < \lambda$ are such that:
\medskip
\roster
\item "$(*)(1)$"  $PT(\kappa,\aleph_0)$, $\kappa$ regular
\item "(2)"  $NPT(\mu,\aleph_0)$, $(\mu \text{ regular })$
\item "(3)"  $PT(\lambda,\aleph_0)$, $\lambda = \lambda^{< \lambda}$, even
if we add a Cohen subset to $\lambda$
\item "(4)"  $\lambda \in I[\lambda]$ (see on it and references in the
discussion after \scite{2.2}; if $\lambda$ is strongly inaccessible or
$\lambda = \chi^+ \and \chi = \chi^{< \chi}$ this always holds).
\endroster
\medskip

\noindent
(So all three are regular uncountable cardinals).
\newline
Shoot a non-reflecting stationary subset ${\underset\sim\to S}$
of $\{ \delta < \lambda:\text{cf}(\delta) = \kappa \}$ by the forcing notion

$$
\align
Q_0 = \biggl\{
h:&h \text{ a function from } (\alpha +1) \text{ to } \{ 0,1 \},
\alpha < \lambda,\\
  &[\delta \le \alpha \and h(\delta) = 1 \Rightarrow cf(\delta) = \kappa],\\
  &h^{-1}(\{ 1 \}) \text{ is non-reflecting} \biggr\}.
\endalign
$$
\medskip

\noindent
In $V^{Q_0}$ no bounded subset of $\lambda$ is added (so $NPT(\mu,\aleph_0)$,
$PT(\kappa,\aleph_0)$ continue to hold).
Also trivially $NPT(\lambda,\kappa)$ hence
$NPT(\lambda,\mu)$ (as $\mu \ge \kappa$, e.g. for $\delta \in
{\underset\sim\to S}$ choose $A_\delta \subseteq \delta$,
$|A_\delta| = \text{ otp}(A_\delta) = \kappa$,
$\text{ sup}(A_\delta) = \delta$).

So it is enough to prove $PT(\lambda,\aleph_0)$.  If not, there is a
$\lambda$-set ${\frak S}$ and $\{ A_\eta:\eta \in {\frak S}_f \}$ 
witnessing it as in \cite{Sh:161, {\S3}} or see \cite{Sh:521, {AP}}.  As
$PT(\kappa,\aleph_0)$ (also in $V^Q$) we have $\eta \in S \Rightarrow
\lambda_\eta \ne \kappa$, hence (see \cite{Sh:161, {\S3}} or
\cite{Sh:521, {AP}})
\medskip
\roster
\item "$(*)$"  $<i> \, \in S \Rightarrow cf(i) \ne \kappa$.
\endroster
\medskip

\noindent
So if we kill the stationary set ${\underset\sim\to S}$ by a forcing notion
\newline
${\underset\sim\to Q _1} = \{ g:a \text{ a function for }
\alpha + 1 < \lambda \text{ to } \{ 0,1 \}, h^{-1}(\{ 1 \})
\text{ closed disjoint to } {\underset\sim\to S} \}$, then in the universe
$\left( V^{Q_0} \right)^{\underset\sim\to Q _1}$,
 still $NPT(\lambda,\aleph_0)$ by the same example. \newline
[Why?  The set $S = \{i < \lambda:<i> \in {\frak S}\}$ is stationary
(as ${\frak S}$ is a $\lambda$-set).  Now (see \cite{Sh:88, {AP}}) the set
$S$ is in $I[\lambda]$ in $V^Q$ (as $\lambda$ is) and is stationary and
disjoint to ${\underset\sim {}\to S_0}$ hence is stationary also in
$\left( V^{Q_0} \right)^{\underset\sim {}\to Q_1}$, and no bounded subset of
$\lambda$ is added.]

But $Q_0*{\underset\sim\to Q_1}$ is just adding a Cohen subset of $\lambda$
hence by $(*)(3)$ we know 
$\left( V^{Q_0} \right)^{*{\underset\sim {}\to Q_1}} \models PT(\lambda,
\aleph_0)$.
So we have gotten an example as required in the lemma, provided we start with
a universe satisfying $(*)$.  So when $(*)$ holds?
\enddemo
\bigskip

\demo{Case 1}  $\lambda$ supercompact, $\kappa$ first measurable, choose 
any $\mu$ such that $\kappa < \mu < \lambda$, \newline
$\mu = \text{ cf}(\mu)$, we 
can add to $\mu$ a set $S \subseteq \{ \delta < \mu:\text{cf}(\delta) 
= \aleph_0 \}$ stationary non-reflecting and then make $\lambda$ Laver
indestructible by forcing notions not adding subsets of $\lambda$. \newline
So G.C.H. is O.K.
\enddemo
\bigskip

\demo{Case 2 (Small Cardinals)}  By the forcing in \cite{MgSh:204} (but for
$\lambda$ regular adding $\lambda$-Cohen sets) the following case is O.K.

$$
GCH + \kappa = \aleph_{\omega^2+1} + \mu = \aleph_{\omega_1+1} +
\lambda = \aleph_{\omega_1+\omega^2+1}.
$$
\enddemo
\bigskip

\demo{Case 3}  $V = L,\kappa < \lambda$ are uncountable weakly compact 
cardinals in $L,\mu = \kappa^+,P$ is a $\mu^+$-complete forcing
notion making $(*)(3)$ hold without adding any subset to $\mu$.
(E.g. 
$P$ is $P_\lambda$, where ${\langle P_i,{\underset\sim \to Q_j}:i \le \lambda,
j < \lambda \rangle}$ is an iteration with Easton
support, $\lambda_j$-the jth inaccessible cardinal in $L$ which is
$> \kappa$, $Q_j$ is the $\lambda_j$-Cohen forcing in $V^{P_j}$); 
i.e. \newline
$Q_j = \{ h:h \text{ a function from some } \alpha < \lambda_j \text{ to }
\lambda_j \text{ in } V^{P_j} \}$. \newline
Now in $V^P, \kappa, \mu, \lambda$ are as required. 
\hfill$\square_{\scite{3.1}}$
\enddemo
\bigskip

\noindent
By work of Jensen there is a non-reflecting stationary 
$S \subseteq   \lambda$  for  $\lambda$ regular non weakly compact

\newpage

REFERENCES.  

\item{[EM]} P.~Eklof and A.~Mekler.
\relax {\it {Almost free modules; Set theoretic methods}}.
\relax North Holland Library, 1990.

\item{[MgSh 204]} Menachem Magidor and Saharon Shelah.
\relax {When does almost free imply free? (For groups, transversal etc.)}.
\relax {\it {Journal of the American Mathematical Society}}, {\bf accepted}.

\item{[Sh 521]} Saharon Shelah.
\relax {If there is an exactly $\lambda$-free abelian group then there is a
  $\lambda$-separable one}.
\relax {\it Journal of Symbolic Logic}, {\bf to appear}.

\item{[Sh 460]} Saharon Shelah.
\relax {The Generalized Continuum Hypothesis revisited}.
\relax {\it {Israel Journal of Mathematics}}, {\bf submitted}.

\item{[Sh 108]} Saharon Shelah.
\relax {On successors of singular cardinals}.
\relax In {\it {Logic Colloquium '78 (Mons, 1978)}}, volume~97 of {\it
  {Stud. Logic Foundations Math}}, pages 357--380. {North-Holland,
  Amsterdam-New York}, 1979.

\item{[Sh 161]} Saharon Shelah.
\relax {Incompactness in regular cardinals}.
\relax {\it {Notre Dame Journal of Formal Logic}}, {\bf 26}:195--228, 1985.

\item{[Sh 88]} Saharon Shelah.
\relax {Classification of nonelementary classes. II. Abstract elementary
  classes}.
\relax In {\it {Classification theory (Chicago, IL, 1985)}}, volume 1292 of
  {\it {Lecture Notes in Mathematics}}, pages 419--497. {Springer, Berlin},
  1987.
\relax {Proceedings of the USA--Israel Conference on Classification Theory,
  Chicago, December 1985; ed. Baldwin, J.T.}

\item{[Sh 345]} Saharon Shelah.
\relax {Products of regular cardinals and cardinal invariants of products of
  Boolean algebras}.
\relax {\it {Israel Journal of Mathematics}}, {\bf 70}:129--187, 1990.

\item{[Sh 420]} Saharon Shelah.
\relax {Advances in Cardinal Arithmetic}.
\relax In {\it {Finite and Infinite Combinatorics in Sets and Logic}}, pages
  355--383. Kluwer Academic Publishers, 1993.
\relax {N.W. Sauer et al (eds.)}.

\item{[Sh 355]} Saharon Shelah.
\relax {$\aleph _{\omega +1}$ has a Jonsson Algebra}.
\relax In {\it {Cardinal Arithmetic}}, volume~29 of {\it {Oxford Logic
  Guides}}, chapter~II. {Oxford University Press}, 1994.

\item{[Sh:g]} Saharon Shelah.
\relax {\it {Cardinal Arithmetic}}, volume~29 of {\it {Oxford Logic
  Guides}}.
\relax {Oxford University Press}, 1994.

\end{thebibliography}
\shlhetal
\enddocument
\bye

%% file: mathdefs.tex
\expandafter\ifx\csname mathdefs.tex\endcsname\relax
  \expandafter\gdef\csname mathdefs.tex\endcsname{}
\else \message{Hey!  Apparently you were trying to
  \string\input{mathdefs.tex} twice.   This does not make sense.} 
\errmessage{Please edit your file (probably \jobname.tex) and remove
any duplicate ``\string\input'' lines} \fi




\catcode`\X=12\catcode`\@=11

\def\n@wcount{\alloc@0\count\countdef\insc@unt}
\def\n@wwrite{\alloc@7\write\chardef\sixt@@n}
\def\n@wread{\alloc@6\read\chardef\sixt@@n}
\def\r@s@t{\relax}\def\v@idline{\par}\def\@mputate#1/{#1}
\def\l@c@l#1X{\firstpart.#1}\def\gl@b@l#1X{#1}\def\t@d@l#1X{{}}

\def\crossrefs#1{\ifx\all#1\let\tr@ce=\all\else\def\tr@ce{#1,}\fi
   \n@wwrite\cit@tionsout\openout\cit@tionsout=\jobname.cit 
   \write\cit@tionsout{\tr@ce}\expandafter\setfl@gs\tr@ce,}
\def\setfl@gs#1,{\def\@{#1}\ifx\@\empty\let\next=\relax
   \else\let\next=\setfl@gs\expandafter\xdef
   \csname#1tr@cetrue\endcsname{}\fi\next}
\def\m@ketag#1#2{\expandafter\n@wcount\csname#2tagno\endcsname
     \csname#2tagno\endcsname=0\let\tail=\all\xdef\all{\tail#2,}
   \ifx#1\l@c@l\let\tail=\r@s@t\xdef\r@s@t{\csname#2tagno\endcsname=0\tail}\fi
   \expandafter\gdef\csname#2cite\endcsname##1{\expandafter
     \ifx\csname#2tag##1\endcsname\relax?\else\csname#2tag##1\endcsname\fi
     \expandafter\ifx\csname#2tr@cetrue\endcsname\relax\else
     \write\cit@tionsout{#2tag ##1 cited on page \folio.}\fi}
   \expandafter\gdef\csname#2page\endcsname##1{\expandafter
     \ifx\csname#2page##1\endcsname\relax?\else\csname#2page##1\endcsname\fi
     \expandafter\ifx\csname#2tr@cetrue\endcsname\relax\else
     \write\cit@tionsout{#2tag ##1 cited on page \folio.}\fi}
   \expandafter\gdef\csname#2tag\endcsname##1{\expandafter
      \ifx\csname#2check##1\endcsname\relax
      \expandafter\xdef\csname#2check##1\endcsname{}%
      \else\immediate\write16{Warning: #2tag ##1 used more than once.}\fi
      \multit@g{#1}{#2}##1/X%
      \write\t@gsout{#2tag ##1 assigned number \csname#2tag##1\endcsname\space
      on page \number\count0.}%
   \csname#2tag##1\endcsname}}
\def\multit@g#1#2#3/#4X{\def\t@mp{#4}\ifx\t@mp\empty%
      \global\advance\csname#2tagno\endcsname by 1 
      \expandafter\xdef\csname#2tag#3\endcsname
      {#1\number\csname#2tagno\endcsnameX}%
   \else\expandafter\ifx\csname#2last#3\endcsname\relax
      \expandafter\n@wcount\csname#2last#3\endcsname
      \global\advance\csname#2tagno\endcsname by 1 
      \expandafter\xdef\csname#2tag#3\endcsname
      {#1\number\csname#2tagno\endcsnameX}
      \write\t@gsout{#2tag #3 assigned number \csname#2tag#3\endcsname\space
      on page \number\count0.}\fi
   \global\advance\csname#2last#3\endcsname by 1
   \def\t@mp{\expandafter\xdef\csname#2tag#3/}%
   \expandafter\t@mp\@mputate#4\endcsname
   {\csname#2tag#3\endcsname\lastpart{\csname#2last#3\endcsname}}\fi}
\def\t@gs#1{\def\all{}\m@ketag#1e\m@ketag#1s\m@ketag\t@d@l p
   \m@ketag\gl@b@l r \n@wread\t@gsin
   \openin\t@gsin=\jobname.tgs \re@der \closein\t@gsin
   \n@wwrite\t@gsout\openout\t@gsout=\jobname.tgs }
\outer\def\localtags{\t@gs\l@c@l}
\outer\def\globaltags{\t@gs\gl@b@l}
\outer\def\newlocaltag#1{\m@ketag\l@c@l{#1}}
\outer\def\newglobaltag#1{\m@ketag\gl@b@l{#1}}

\newif\ifpr@ 
\def\m@kecs #1tag #2 assigned number #3 on page #4.%
   {\expandafter\gdef\csname#1tag#2\endcsname{#3}
   \expandafter\gdef\csname#1page#2\endcsname{#4}
   \ifpr@\expandafter\xdef\csname#1check#2\endcsname{}\fi}
\def\re@der{\ifeof\t@gsin\let\next=\relax\else
   \read\t@gsin to\t@gline\ifx\t@gline\v@idline\else
   \expandafter\m@kecs \t@gline\fi\let \next=\re@der\fi\next}
\def\pretags#1{\pr@true\pret@gs#1,,}
\def\pret@gs#1,{\def\@{#1}\ifx\@\empty\let\n@xtfile=\relax
   \else\let\n@xtfile=\pret@gs \openin\t@gsin=#1.tgs \message{#1} \re@der 
   \closein\t@gsin\fi \n@xtfile}

\newcount\sectno\sectno=0\newcount\subsectno\subsectno=0
\newif\ifultr@local \def\ultralocal{\ultr@localtrue}
\def\firstpart{\number\sectno}
\def\lastpart#1{\ifcase#1 \or a\or b\or c\or d\or e\or f\or g\or h\or 
   i\or k\or l\or m\or n\or o\or p\or q\or r\or s\or t\or u\or v\or w\or 
   x\or y\or z \fi}

\def\resetall{\global\advance\sectno by 1\subsectno=0
   \gdef\firstpart{\number\sectno}\r@s@t}
\def\resetsub{\global\advance\subsectno by 1
   \gdef\firstpart{\number\sectno.\number\subsectno}\r@s@t}
\def\newsection#1\par{\resetall\vskip0pt plus.3\vsize\penalty-250
   \vskip0pt plus-.3\vsize\bigskip\bigskip
   \message{#1}\leftline{\bf#1}\nobreak\bigskip}
\def\subsection#1\par{\ifultr@local\resetsub\fi
   \vskip0pt plus.2\vsize\penalty-250\vskip0pt plus-.2\vsize
   \bigskip\smallskip\message{#1}\leftline{\bf#1}\nobreak\medskip}

\def\t@gsoff#1,{\def\@{#1}\ifx\@\empty\let\next=\relax\else\let\next=\t@gsoff
   \def\@@{p}\ifx\@\@@\else
   \expandafter\gdef\csname#1cite\endcsname##1{\zeigen{##1}}
   \expandafter\gdef\csname#1page\endcsname##1{?}
   \expandafter\gdef\csname#1tag\endcsname##1{\zeigen{##1}}\fi\fi\next}
\def\verbatimtags{\ifx\all\relax\else\expandafter\t@gsoff\all,\fi}
\def\zeigen#1{\hbox{$\langle$}#1\hbox{$\rangle$}}

\def\(#1){\edef\dot@g{\ifmmode\ifinner(\hbox{\noexpand\etag{#1}})
   \else\noexpand\eqno(\hbox{\noexpand\etag{#1}})\fi
   \else(\noexpand\ecite{#1})\fi}\dot@g}

\newif\ifbr@ck
\def\eat#1{}
\def\[#1]{\br@cktrue[\br@cket#1'X]}
\def\br@cket#1'#2X{\def\temp{#2}\ifx\temp\empty\let\next\eat
   \else\let\next\br@cket\fi
   \ifbr@ck\br@ckfalse\br@ck@t#1,X\else\br@cktrue#1\fi\next#2X}
\def\br@ck@t#1,#2X{\def\temp{#2}\ifx\temp\empty\let\neext\eat
   \else\let\neext\br@ck@t\def\temp{,}\fi
   \def\teemp{#1}\ifx\teemp\empty\else\rcite{#1}\fi\temp\neext#2X}
\def\resetbr@cket{\gdef\[##1]{[\rtag{##1}]}}
\def\references{\resetbr@cket\newsection References\par}

\newtoks\symb@ls\newtoks\s@mb@ls\newtoks\p@gelist\n@wcount\ftn@mber
    \ftn@mber=1\newif\ifftn@mbers\ftn@mbersfalse\newif\ifbyp@ge\byp@gefalse
\def\defm@rk{\ifftn@mbers\n@mberm@rk\else\symb@lm@rk\fi}
\def\n@mberm@rk{\xdef\m@rk{{\the\ftn@mber}}%
    \global\advance\ftn@mber by 1 }
\def\rot@te#1{\let\temp=#1\global#1=\expandafter\r@t@te\the\temp,X}
\def\r@t@te#1,#2X{{#2#1}\xdef\m@rk{{#1}}}
\def\b@@st#1{{$^{#1}$}}\def\str@p#1{#1}
\def\symb@lm@rk{\ifbyp@ge\rot@te\p@gelist\ifnum\expandafter\str@p\m@rk=1 
    \s@mb@ls=\symb@ls\fi\write\f@nsout{\number\count0}\fi \rot@te\s@mb@ls}
\def\byp@ge{\byp@getrue\n@wwrite\f@nsin\openin\f@nsin=\jobname.fns 
    \n@wcount\currentp@ge\currentp@ge=0\p@gelist={0}
    \re@dfns\closein\f@nsin\rot@te\p@gelist
    \n@wread\f@nsout\openout\f@nsout=\jobname.fns }
\def\m@kelist#1X#2{{#1,#2}}
\def\re@dfns{\ifeof\f@nsin\let\next=\relax\else\read\f@nsin to \f@nline
    \ifx\f@nline\v@idline\else\let\t@mplist=\p@gelist
    \ifnum\currentp@ge=\f@nline
    \global\p@gelist=\expandafter\m@kelist\the\t@mplistX0
    \else\currentp@ge=\f@nline
    \global\p@gelist=\expandafter\m@kelist\the\t@mplistX1\fi\fi
    \let\next=\re@dfns\fi\next}
\def\symbols#1{\symb@ls={#1}\s@mb@ls=\symb@ls} 
\def\bigsymbol{\textstyle}
\symbols{\bigsymbol\ast,\dagger,\ddagger,\sharp,\flat,\natural,\star}
\def\ftnumbers{\ftn@mberstrue} \def\ftsymbols{\ftn@mbersfalse}
\def\paginal{\byp@ge} \def\resetftnumbers{\ftn@mber=1}
\def\ftnote#1{\defm@rk\expandafter\expandafter\expandafter\footnote
    \expandafter\b@@st\m@rk{#1}}

\long\def\jump#1\endjump{}
\def\ssum{\mathop{\lower .1em\hbox{$\textstyle\Sigma$}}\nolimits}

\def\qed{\nobreak\kern 1em \vrule height .5em width .5em depth 0em}
\def\newneq{\hbox{\rlap{\hbox to 1\wd9{\hss$=$\hss}}\raise .1em 
   \hbox to 1\wd9{\hss$\scriptscriptstyle/$\hss}}}
\def\subsetne{\setbox9 = \hbox{$\subset$}\mathrel{\hbox{\rlap
   {\lower .4em \newneq}\raise .13em \hbox{$\subset$}}}}
\def\supsetne{\setbox9 = \hbox{$\subset$}\mathrel{\hbox{\rlap
   {\lower .4em \newneq}\raise .13em \hbox{$\supset$}}}}

\def\vbar{\mathchoice{\vrule height6.3ptdepth-.5ptwidth.8pt\kern-.8pt}
   {\vrule height6.3ptdepth-.5ptwidth.8pt\kern-.8pt}
   {\vrule height4.1ptdepth-.35ptwidth.6pt\kern-.6pt}
   {\vrule height3.1ptdepth-.25ptwidth.5pt\kern-.5pt}}
\def\f@dge{\mathchoice{}{}{\mkern.5mu}{\mkern.8mu}}
\def\b@c#1#2{{\rm \mkern#2mu\vbar\mkern-#2mu#1}}
\def\b@b#1{{\rm I\mkern-3.5mu #1}}
\def\b@a#1#2{{\rm #1\mkern-#2mu\f@dge #1}}
\def\bb#1{{\count4=`#1 \advance\count4by-64 \ifcase\count4\or\b@a A{11.5}\or
   \b@b B\or\b@c C{5}\or\b@b D\or\b@b E\or\b@b F \or\b@c G{5}\or\b@b H\or
   \b@b I\or\b@c J{3}\or\b@b K\or\b@b L \or\b@b M\or\b@b N\or\b@c O{5} \or
   \b@b P\or\b@c Q{5}\or\b@b R\or\b@a S{8}\or\b@a T{10.5}\or\b@c U{5}\or
   \b@a V{12}\or\b@a W{16.5}\or\b@a X{11}\or\b@a Y{11.7}\or\b@a Z{7.5}\fi}}

\catcode`\X=11 \catcode`\@=12